\documentclass[12pt]{amsart}
\usepackage{amssymb,amsmath,stmaryrd}
\usepackage{hyperref}
\usepackage{mathrsfs}
\hypersetup{pdfpagemode=FullScreen,colorlinks=true}
\def\R{\mathbb{R}}

\def\d{\mathrm{d}_c}
\newtheorem{thm}{Theorem}
\newtheorem{pro}{Proposition}

\newtheorem{lem}[pro]{Lemma}

\newtheorem{exa}[pro]{Example}
\newtheorem{dfi}[pro]{Definition}

\newenvironment{pf}{\begin{trivlist}\item[]{\bf Proof\ }}
{\mbox{}\hfill\rule{.08in}{.08in}\end{trivlist}}

\title[Normal currents and charges on Carnot groups]{Flat compactness of normal currents and charges on Carnot groups}
\author{Antoine Julia}
\thanks{
A.~Julia is supported by the Simons Foundation grant 601941, GD}
\author{Pierre Pansu}
\thanks{P.~Pansu is supported by Agence Nationale de la Recherche, ANR-22-CE40-0004 GOFR}

\address{Antoine Julia:
Universit\'e Paris-Saclay, CNRS, Laboratoire de Ma\-th\'ematiques d'Orsay\\ 91405 Orsay C\'edex, France}
\email{antoine.julia@universite-paris-saclay.fr}

\address{Pierre Pansu:
Universit\'e Paris-Saclay, CNRS, Laboratoire de Ma\-th\'ematiques d'Orsay\\ 91405 Orsay C\'edex, France}
\email{pierre.pansu@universite-paris-saclay.fr}

\begin{document}

\maketitle

\begin{abstract}
We prove that the family of normal currents in the sense of Rumin in a Carnot group is compact in the flat topology. This result is obtained through a dual compactness argument for Rumin forms, using the pseudo-differential calculus in groups developed by Folland, Christ-Geller-G\l{}owacki-Polin and Rumin. As an application, imitating de Pauw-Moonens-Pfeffer, we describe the space of charges on a Carnot group.
\end{abstract}
  
\tableofcontents

\section{Introduction}

In an effort to develop geometric measure theory in Carnot groups, we study an adapted notion of currents, defined by duality with Rumin differential forms. This notion is not new, it plays a crucial role in recent work of D. Vittone on intrinsic Lipschitz graphs in Heisenberg groups, \cite{Vittone2022Sigma}. Our remote goal is to find an appropriate notion and establish properties of integral currents. We are still very far from this goal, and begin modestly with normal currents: we prove that the family of normal currents in the sense of Rumin, with support in a fixed compact set in a Carnot group, is compact in the flat topology. 

In the Euclidean setting, flat compactness of normal currents is usually stated as a by-product of the Deformation Theorem. This theorem (\cite[Theorem 4.2.9]{FedererGMT}) is designed for proving the flat compactness of integral currents. It is not available in a general Carnot group. Instead, we argue by duality. The dual compactness argument for Rumin forms amounts to inverting Rumin's differentials. We rely on M. Rumin's observation that the Rumin complex is maximally hypoelliptic, \cite{Rumin2001Grenoble}, and on the pseudo-differential calculus in groups developed by G. Folland, \cite{Folland1975Subelliptic} and M. Christ, D. Geller, P. G\l owacki and L. Polin, \cite{CGGP1992}. 

As an application of our compactness theorem, imitating T. De Pauw, L. Moonens and W. Pfeffer, \cite{dePauwMoonensPfeffer2009Charges}, we describe the space of charges on a Carnot group. Charges are linear functionals on compactly supported currents, see a precise definition and statement below. 

\subsection{Rumin forms}

A Carnot group is a Lie group $\mathbb{G}$ equipped with a $1$-parameter group $(\delta_t)_{t>0}$ of automorphisms, called \emph{dilations}, with the following property: the Lie algebra $\mathfrak{g}$ is generated at level $1$, i.e. by the $1$-eigenspace of the infinitesimal generator (which is a derivation of $\mathfrak{g}$). This produces a grading $\mathfrak{g}=\bigoplus_{i=1}^s \mathfrak{g}_i$, hence a grading of the exterior algebra $\Lambda^\cdot\mathfrak{g}^*$. One can therefore define weights of covectors, and by extension, of differential forms on $\mathbb{G}$. For instance, a left-invariant form $\omega$ has weight $w$ if and only if $\delta_t^*\omega=t^w \omega$.

\begin{exa}\textbf{Differential forms on the $3$-dimensional Heisenberg group}. There, $0$-forms have weight $0$, $1$-forms exist in weights $1$ and $2$, $1$-forms of weight $2$ being those which vanish along the left-invariant contact structure. Similarly, $2$-forms exist in weights $2$ and $3$, $2$-forms of weight $2$ being generated by the differential of the contact form, $2$-forms of weight $3$ being those which vanish along the left-invariant contact structure. $3$-forms have weight $4$.
\end{exa}

Rumin's theory aims at improving de Rham's complex on a Carnot group to make it compatible with dilations. Ideally, by reducing it to exactly one weight per degree, like in abelian groups. There is an obstruction, the grading of the cohomology of left-invariant forms. Rumin constructs an optimal complex, homotopy equivalent to de Rham's complex, with a minimal number of weights in each degree (equal to the number of weights in the cohomology). This complex can be viewed as a subspace of differential forms, defined by the vanishing of certain components, equipped with a modified differential $\d$. $\d$ is a differential operator whose order varies with the degree and weight of forms. It satisfies $\d\circ \d=0$. Rumin's complex is homotopy equivalent to de Rham's complex, so it can be used to compute cohomology.

\begin{exa}\textbf{Rumin forms on the $3$-dimensional Heisenberg group}. Rumin's construction selects $0$-forms of weight $0$, $1$-forms of weight $1$, $2$-forms of weight $3$ and $3$-forms of weight $4$. Rumin's differential $\d$ has order $2$ in degree $1$ and $1$ in other degrees. On $0$-forms, $\d$ is the restriction of the usual differential to the left-invariant contact structure. 
\end{exa}

Heisenberg groups are examples where the procedure is fully successful, with one weight in each degree. Here is a less successful (but still optimal) example.

\begin{exa}\textbf{Rumin forms on the $4$-dimensional Engel group}. Rumin's construction selects $0$-forms of weight $0$, $1$-forms of weight $1$, $2$-forms of weights $3$ and $4$, $3$-forms of weight $6$ and $4$-forms of weight $7$. Rumin's differential $\d$ has two components of respective orders $2$ and $3$ in degrees $1$ and again in degree $2$. It has order $1$ in other degrees.
\end{exa}

Rumin's construction is not fully invariant: a choice of left-invariant Riemannian metric on the Carnot group, adapted to the eigenspaces of the dilations, is needed. This choice also allows to normalize Haar measure and to measure the pointwise norm of Rumin forms. Furthermore, it determines a left-invariant distance, the Carnot-Carath\'eodory distance, which is homogeneous of degree $1$ under dilations.

For a survey of the theory, see \cite{Rumin2005Srni}. Slightly different approaches can be found in the following recent sources: \cite{LerarioTripaldi2022}, \cite{FicherTripaldi2022}.

\subsection{Rumin currents}

A \emph{Rumin current} is a continuous functional on the space of smooth compactly supported Rumin forms. The boundary operator $\partial_c$ is defined by duality, 
$$
\langle \partial_c T,\omega \rangle =\langle T,\d\omega \rangle.
$$
The \emph{support} of $T$ is such that $T$ vanishes on forms with support in the complement. The \emph{mass} of a current is defined by duality with the $C^0$ norm on smooth compactly supported Rumin forms,
$$
\mathbf{M}(T)=\sup\{T(\omega)\,;\,\|\omega\|_{C^0(\mathbb{G})}\le 1\}.
$$
A Rumin current $T$ is \emph{normal} if $T$ and $\partial_c T$ have finite mass: $\mathbf{N}(T)=\mathbf{M}(T)+\mathbf{M}(\partial_c T)$.

The \emph{flat mass} of a Rumin current $T$ is the infimum of $\mathbf{M}(R)+\mathbf{M}(S)$ over all expressions $T=S+\partial_c R$.

\subsection{Examples}

\textbf{Diffuse currents}. A $C^1$ Rumin form $\phi$ defines a current $P(\phi)$ of complementary dimension by
$$
\langle P(\phi),\omega \rangle = \int\phi\wedge\omega.
$$
(Recall that Rumin forms are not stable by the exterior product in general, but for forms of complementary degree it is well defined.) 
Then $\mathbf{M}(P(\phi))=\|\phi\|_1$, $\partial_c P(\phi)=\pm P(\d\phi)$. Such diffuse currents $P(\phi)$ are dense in flat norm in the space of normal Rumin currents. 

\bigskip

\textbf{Currents of integration}. In the $3$-dimensional Heisenberg group $\mathbb{H}^1$, a smooth 2-submanifold $V$ with boundary defines a Rumin current $T_V$ of mass equal to its Hausdorff $3$-dimensional measure. The Rumin boundary $\partial_c T_V$ differs from $\partial T_V$: if $\omega$ is a smooth compactly supported $1$-form,
$$
\langle \partial_c T_V,\omega \rangle = \langle \partial T_V,\omega-\frac{d\omega}{d\theta}\theta \rangle,
$$
where $\theta$ denotes the unit vertical left-invariant form. Therefore \break $M(\partial_c T_V)<\infty$ if and only if $\partial V$ is a horizontal curve. In this case, $\partial_c T_V=\partial T_V$.

More generally, in a Carnot group, a compact smooth submanifold with boundary whose boundary is horizontal defines a current $T_V$ which is a normal Rumin current, and $\partial_c T_V=\partial T_V$. (This is, by the way, a reasonable candidate for an integral Rumin current.)

\subsection{Results}\label{sec:intro-results}

\begin{thm}[Compactness]
\label{thm:compactness}
The space of normal Rumin currents with bounded normal mass and support in a fixed compact subset of $\mathbb{G}$ is compact in flat topology.
\end{thm}

As an application, following T. de Pauw, L. Moonens and W. Pfeffer, we define \emph{Rumin charges} on a Carnot group $\mathbb{G}$ as linear functionals on the space of compactly supported normal Rumin currents, which are continuous on each 
$$
S(K,\nu)=\{\text{Rumin currents with support in } K \text{ and normal mass } \le \nu\},
$$
equipped with the flat mass, for all compact sets $K\subset\mathbb{G}$ and $\nu>0$.

The following is a direct generalization of \cite{dePauwMoonensPfeffer2009Charges} to Carnot groups, the only new ingredient needed being Theorem \ref{thm:compactness}.

\begin{thm}[Representation of Rumin charges]
\label{thm:charges}
The space of Rumin charges is $C^0+\d C^0$, i.e. if $\phi$ is a continuous Rumin form on $\mathbb{G}$, both $\phi$ and $\d\phi$ define Rumin charges. Conversely, every Rumin charge is of the form $\phi+\d\psi$, where $\phi$ and $\psi$ are continuous Rumin forms on $\mathbb{G}$.
\end{thm}

\section{Preliminaries}

\subsection{Functional analysis}
\begin{pro}\label{prop:dual-compactness}
Let $X$ and $Y$ be normed vectorspaces. Let $P:X\to Y$ be a compact operator, i.e. bounded sets are mapped to precompact sets. Then $P^*:Y^*\to X^*$ is compact as well. 
\end{pro}

\begin{pf}
Let $B$ denote the unit ball in $Y$ and $B^*$ the unit ball in $Y^*$, i.e. the set of linear functionals $y^*$ on $Y$ such that $\sup_{y\in B} \langle y^*,y \rangle \leq 1$. Every element of $B^*$ is a $1$-Lipschitz function on $Y$ which sends $0$ to $0$. Hence $B^*$ can be viewed as a subset $S$ of the set of $1$-Lipschitz functions on the completion of $P(B)$ which send $0$ to $0$. Since $P(B)$ is a precompact metric space, its completion (which need not coincide with its closure in $Y$) $\overline{P(B)}$ is compact. According to the Arzel\`a-Ascoli Theorem, $S$ is precompact with respect to uniform convergence on $\overline{P(B)}$. Note that, for $y^*\in Y^*$,
\begin{align*}
\|y^*\|_{C^0(\overline{P(B)})}
&=\|y^*\|_{C^0(P(B))}
=\sup_{x\in B}\langle y^*,P(x) \rangle
=\sup_{x\in B}\langle P^*(y^*),x \rangle\\
&=\|P^*(y^*)\|_{X^*}.
\end{align*}
Hence $P^*(B^*)$ is precompact in the norm of $X^*$.
\end{pf}

\subsection{Rumin forms and currents}

The space of smooth Rumin forms of degree $m$ with compact support is denoted by $\mathscr{D}_c^m(\mathbb{G})$, endowing it with the $C^\infty$ topology yields a Fr\'echet space. We will be using three different semi-norms on this space: the $C^0$ norm $\omega \mapsto |\cdot|$, the flat norm 
\[ 
\omega \mapsto \mathbf{F}(\omega):=\max \{|\omega|,  |\d \omega|\}
\]
and the ``normal norm''
\begin{align*}
\mathbf{N}(\omega):= \inf\{ \max\{ |\phi| , |\psi|\}\,;\, \omega= \phi+ \d \psi, \phi \in \mathscr{D}_c^m(\mathbb{G}), \psi \in \mathscr{D}_c^{m-1}(\mathbb{G})\}.
\end{align*}

The dual space of $\mathscr{D}_c^m(\mathbb{G})$ for the Fr\'echet topology is the space of \emph{Rumin currents} of dimension $m$, $\mathscr{D}_{c,m}(\mathbb{G})$. The boundary of a Rumin current $T$ of dimension $m\ge 1$ is the Rumin current of dimension $m-1$ defined by $\partial_c T (\omega) = T(\d \omega)$. To a current $T\in \mathscr{D}_{c,m}(\mathbb{G})$, we can associate three ``norms'' (which are not always finite),
\begin{enumerate}
  \item The \emph{mass} 
  \[
     \mathbf{M}(T):= \sup \{T(\omega), |\omega|\le 1\}.
  \]
  \item The \emph{normal mass} 
  \[
        \mathbf{N}(T) :=\mathbf{M}(T) + \mathbf{M}(\partial T) =  \sup \{T(\omega), \mathbf{N}(\omega)\le 1\}.
  \]
  \item The \emph{flat mass}
  \[
    \mathbf{F}(T) := \sup \{T(\omega), \mathbf{F}(\omega)\le 1\} = \inf\{\mathbf{M}(R) + \mathbf{M}(S), T= S+\partial_c R\}.
    \]
  \end{enumerate}
    The identity in the third definition holds provided the flat mass of $T$ is finite.
   Let us prove the identities in the second and third definition. For the normal mass, the direction $\ge$ is straightforward, for the opposite direction we can suppose that $\mathbf{M}(T) + \mathbf{M}(\partial_c T) < +\infty$. In particular there are forms $\phi$ and $\psi$ with $|\phi|=|\psi|=1$ such that $T(\phi)$ and $\partial T(\psi)$ are arbitrarily close to $\mathbf{M}(T)$ and $\mathbf{M}(\partial T)$ respectively. It suffices to then take $\omega:=\phi + \d \psi$. For the flat mass, inequality $\le$ is clear. To prove the opposite inequality, one needs to use the Hahn-Banach Theorem as in \cite[4.1.12]{FedererGMT}.
   
    \medskip
    
  Finally, the \emph{support} of a current $T$ is the smallest closed subset $K$ of $\mathbb{G}$ such that $T(\omega)=0$ whenever $\omega$ is supported in a compact subset of the complement of $K$. We denote by  $\mathbf{N}_{c,m}(\mathbb{G})$ the subspace of \emph{normal Rumin currents} of dimension $m$ with finite normal mass and compact support and by $\mathbf{F}_{c,m}(\mathbb{G})$, the closure of $\mathbf{N}_{c,m}(\mathbb{G})$ with respect to the flat mass.

  \subsection{Continuity of embeddings}
  
  Our aim is to prove that the space of normal Rumin currents of dimension $m$ supported in a compact set $K\subset \mathbb{G}$, $\mathbf{N}_{c,m, K}(\mathbb{G})$, embeds compactly in $\mathbf{F}_{c,m,K}(\mathbb{G})$ for the flat mass topology. In order to do this, we first prove that the following chain of maps
  \begin{equation}\label{eq:embeddingchain}
      \mathbf{N}_{c,m,K}(\mathbb{G}) \hookrightarrow E_2^* \stackrel{\Theta^*}{\hookrightarrow} E_1^* \hookrightarrow \mathbf{F}_{c,m}(\mathbb{G}),
    \end{equation}
is continuous where $E_2^*$ and $E_1^*$ are dual spaces to two Banach spaces of differential forms such that the map $\Theta : E_1\hookrightarrow E_2$ is compact (this last fact 'will be proved in the following section). By Proposition~\ref{prop:dual-compactness}, the dual arrow is compact and thus the composition of the three arrows in \eqref{eq:embeddingchain} is compact as well. Let us now define these spaces. 

\medskip

  Let $B'$ be an open ball in $\mathbb{G}$  containing a neighbourhood of $K$. Let $B$ be a larger concentric open ball. Let $\mathscr{E}^m_{c}(B)$ denote the space of smooth Rumin forms on $B$ which are bounded on $B$ and whose Rumin differential is bounded on $B$. On this space, let 
  $$
\mathbf{F}_B=\max (|\cdot |_{C^0(B)}, |\d \cdot|_{C^0(B})),\quad \mathbf{F}_K=\max (|\cdot |_{C^0(K)}, |\d \cdot|_{C^0(K)}).
  $$
  
  The space $E_1$ is the completion of the quotient space
    \[
    \big( \mathscr{E}^m_{c}(B),\mathbf{F}_B \big )/ \big\{\omega \in  \mathscr{E}^m_{c}(B), \mathbf{F}_K(\omega) =0 \big \}.
  \]

  To define $E_2$, consider the flat semi-norms on $\mathscr{E}^m_{c}(B')$ defined by
\begin{align*}
\mathbf{N}_{B'}(\omega) &= \inf \big \{ \max (  |\phi|_{C^0(B')},|\psi|_{C^0(B')})\,;\,\phi\in \mathscr{E}^m_{c}(B'),\,\psi\in \mathscr{D}^{m-1}_{c,B'},\\
& \hskip1cm \omega = \phi + \d \psi \big \}.\\
 \mathbf{N}_{K}(\omega) &= \inf \big \{ \max (  |\phi|_{C^0(K)} , |\psi|_{C^0(K)})\,;\,\phi\in \mathscr{E}^m_{c}(B'),\,\psi\in \mathscr{D}^{m-1}_{c,B'},\\
 & \hskip1cm  \omega = \phi + \d \psi \big \}.
\end{align*}
   The space $E_2$ is then the completion of the quotient 
    \begin{align*}
       \big (\mathscr{E}^m_{c}(B'), \mathbf{N}_{B'} \big ) / \{\omega\in \mathscr{E}^m_{c}(B')\,:\, \mathbf{N}_{K}(\omega) = 0 \}.
    \end{align*}
    Recall that the canonical norm on a quotient of a normed space by a closed subspace is the infimum of the norms of its representatives.
    
    The map $\Theta:E_1\to E_2$ is induced by restricting forms defined on the larger ball $B$ to the smaller ball $B'$.
    
  \begin{pro}
    The spaces $E_1$ and $E_2$ are Banach spaces. Furthermore, the arrows
    $$
\mathbf{N}_{c,m,K}(\mathbb{G}) \hookrightarrow E_2^* \quad\text{and}\quad E_1^* \hookrightarrow \mathbf{F}_{c,m}(\mathbb{G})
    $$
    are continuous.
  \end{pro}
  
  \begin{pf}
  
    Since $\mathbf{F}_K \le \mathbf{F}_B$ and $\mathbf{N}_K \le \mathbf{N}_{B'}$, the subspaces by which we quotient are closed, so $E_1$ and $E_2$ are Banach spaces. 
    
\medskip

Let us check the continuity of the first arrow. The key point is that if $T$ is a normal current with support in $K$, then $T$ and $\partial_c T$ are representable by integration and carried by measures supported in $K$, thus
$$
T(\phi)\le \mathbf{M}(T)|\phi|_{C^0(K)},\quad \partial_c T(\phi)\le \mathbf{M}(\partial_c T)|\phi|_{C^0(K)}
$$
for all smooth forms defined in a neighborhood of $K$. It first implies that $T$ vanishes on $\{\omega\in \mathscr{E}^m_{c}(B')\,;\, \mathbf{N}_{K}(\omega) = 0 \}$, and so defines an element of $E_2^*$. Second, that when $\omega=\phi+\d\psi$ with $\phi,\psi$ defined on $B'$, 
    \begin{align*}
T(\omega)&\le \mathbf{M}(T)|\phi|_{C^0(K)}+\mathbf{M}(\partial_c T)|\psi|_{C^0(K)}\\
&\le \mathbf{N}(T)\max(|\phi|_{C^0(B')},|\psi|_{C^0(B')}),
    \end{align*}
which leads to $\|T\|_{E_2^*}\le \mathbf{N}(T)$.

\medskip

Finally, let us check the continuity of the second arrow. Smooth compactly supported Rumin forms restrict to elements of $\mathscr{E}^m_{c}(B)$, $\omega\mapsto \omega_{|B}$. Obviously,
$$
\mathbf{F}_{B}(\omega_{|B})\le \mathbf{F}(\omega),
$$
so elements of $E_1^*$ define currents of flat norm less than their $E_1^*$-norm.

\end{pf}

  \section{Compactness of forms}
  
  The aim of this section is to prove the following result.
  
  \begin{thm}\label{thm:forms}
    The map $\Theta:E_1 \hookrightarrow E_2$ induced by restricting forms defined on $B$ to $B'$ is compact.
  \end{thm}
  
  This result will follow from the construction of a partial inverse to the Rumin differential, together with some control on its norm.  We need to introduce some results of pseudo-differential calculus on Carnot groups, adapted to the Rumin complex. 
  
  \subsection{Global inversion of the Rumin complex}
  
Following \cite{CGGP1992} we define an operator on the space $\mathcal{S}_0$ of forms with coefficients in the Schwartz class, all of whose polynomial moments vanish. This operator will be extended to Sobolev spaces of forms later on. Following \cite{Folland1975Subelliptic} and \cite{PansuRumin2018AHL}, we also consider the square root $|\nabla_H|$ of the scalar subLaplacian $\Delta_H$ on functions on $\mathbb{G}$. For each $m\in \R$, the pseudo-differential operator $|\nabla_H|^m$ is homogeneous of order $m$. We fix once and for all a left-invariant frame, and let $|\nabla_H|^m$ act componentwise on Rumin forms.

    \begin{pro}[{\cite[Proposition 5.2]{PansuRumin2018AHL}}]\label{pro:inverse}
    There exists an operator $K_c$ on the space of Rumin forms in $\mathcal{S}_0$ 
    such that if $\omega$ is such a form,
    \[
      \omega = \d K_c \omega + K_c \d \omega.
    \]
    Furthermore, if $\omega$ is of degree $k$ and pure weight $N$, then $K_c\omega$ is a form of degree $k-1$ and decomposes into forms $(K_c\omega)_\ell$ of pure weight $N-\ell$ for $\ell$ in $\{1,\dots, \min(\delta,N)\}$, where $\delta$ is the largest weight difference between a form $\alpha$ and its Rumin differential $\d\alpha$.
The homogeneous component $(K_c\omega)_\ell$ satisfies for every $m\in \R$ and $p\in (1,+\infty)$,
    \begin{equation}\label{eq:control-lp}
       \big\Vert |\nabla_H|^{m+\ell} (K_c \omega)_\ell\big \Vert_{p} \lesssim \big \Vert |\nabla_H|^{m} \omega \big \Vert_p.
     \end{equation}
   \end{pro}
   
    For positive $m$, the seminorms $\Vert\, |\nabla|^m \cdot \Vert_p$ are comparable to the horizontal Sobolev norms $W_{c}^{1,p}$, which is defined as follows. For a scalar function $u$ on an open set $U\subset\mathbb{G}$, the $W_{c}^{1,p}$ norm is
\begin{align*}
|u|_{W_{c}^{1,p}(U)}=|u|_{L^p(U)}+|\d u|_{L^p(U)},
\end{align*}
and for a form, the $W_{c}^{1,p}$ norm is the sum of the $W_{c}^{1,p}$ norms of its components in the chosen left-invariant frame.
    
    The next statement  is a particular case of \cite[Theorem~4.10]{Folland1975Subelliptic}.

  \begin{pro}\label{pro:sobolev-equivalence}
    Given $p\in (1,\infty)$, for all $f\in \mathcal{S}_0$,
    \begin{equation*}
      \big \Vert f \big \Vert_{W_c^{1,p}} \lesssim \big \Vert |\nabla_H| f \big \Vert_p + \big \Vert f \big \Vert_p.
    \end{equation*}
  \end{pro}
  
    Furthermore, whenever $p$ is larger than the homogeneous dimension $Q$ of $\mathbb{G}$ ($Q=\sum_{i=1}^s i\,\mathrm{dim}(\mathfrak{g}_i)$), the Sobolev $W_c^{1,p}$ norm controls the H\"older semi-norm of order $1-Q/p$ with respect to the Carnot-Carath\'eodory metric. Indeed, following \cite[Theorem~5.15]{Folland1975Subelliptic}, we have
  
  \begin{pro}\label{pro:sobolev-embedding}
    If $f$ is a function or a Rumin form on $\mathbb{G}$, 
    \begin{equation*}
      \sup_{x,y\in \mathbb{G}} \dfrac{|f(y)-f(x)|}{\Vert x^{-1} y \Vert^{1-Q/p}}\lesssim \Vert f \Vert_{W_c^{1,p}}.
    \end{equation*}
  \end{pro}
Let $0<\alpha<1$. Let $U$ be an open subset of $\mathbb{G}$. In the sequel, the space of continuous functions (or Rumin forms) on $U$ whose H\"older semi-norm
  $$
\|f\|_{\dot C_{c}^{\alpha}(U)}:=  \sup_{x,y\in U} \dfrac{|f(y)-f(x)|}{\Vert x^{-1} y \Vert^{\alpha}}
  $$
is finite will be denoted by $\dot C_{c}^{\alpha}(U)$. The H\"older norm is $\|f\|_{ C_{c}^{\alpha}(U)}=\|f\|_{\dot C_{c}^{\alpha}(U)}+\|f\|_{C^{0}(U)}$.

   For negative $m$, we shall need to control the norm of $|\nabla_H|^m f$ with respect to that of $f$. 
   We have the following bound:
   \begin{pro}\label{pro:nablaminus}
     Given $\ell \in (0,Q)$ and $p,q>1$ such that $p^{-1} = q^{-1} - \ell/Q$, for every $f$ in $L^q(\mathbb{G})$, 
     \begin{equation*}
       \Vert |\nabla_H|^{-\ell}f\Vert_p \lesssim \Vert f\Vert_q.
     \end{equation*}
   \end{pro}
   \begin{pf}
By Proposition 3.17 in \cite{Folland1975Subelliptic}, the operator $|\nabla|^{-m}$ has the following kernel
    \[
       R_m(x) := C_m \int_0^\infty t^{(m/2)-1}h(x,t) dt,
     \]
     where $h(x,t)$ is the heat kernel associated to the subLaplacian in $\mathbb{G}$. This means we can write $|\nabla|^{-m} \omega = R_m*\omega$, with $*$ standing for the group convolution in $\mathbb{G}$ (and decomposing $\omega$ in the chosen left-invariant basis). As $m<Q$, by Proposition 3.17 in \cite{Folland1975Subelliptic}, $R_m$ is in $C^\infty (\mathbb{G}\backslash\{0\})$ and is a kernel of type $m$. By Proposition 1.11 in \cite{Folland1975Subelliptic},  given $p,q>1$ such that $p^{-1} = q^{-1}-m/Q$, for $f\in L^q(\mathbb{G})$, 
     $$
     \Vert R_m * f\Vert_p \lesssim_q \Vert f\Vert_q .
     $$ 
   \end{pf}
   
   Putting these estimates together, we get
   
   \begin{lem}\label{lem:kc}
Let $p>Q$. Let $\delta$ be the largest weight increase by the Rumin differential. Then $\delta\le Q-1$. Let $q>1$ be defined by $p^{-1}=q^{-1}-(\delta/Q)$. Then $K_c$ extends to a bounded operator from $L^{p}\cap L^{q}(\mathbb{G})$ to $\dot C_{c}^{1-(Q/p)}(\mathbb{G})\cap L^{p}(\mathbb{G})$.
\end{lem}

\begin{pf}
Combining Propositions \ref{pro:inverse} and \ref{pro:sobolev-equivalence} yields, for a homogeneous differential form $\omega\in\mathcal{S}_0$ of weight $N$ and for every $\ell\in\{1,\ldots,\delta\}$, the homogeneous component of weight $N-\ell$ of $K_c \omega$ satisfies
    \begin{equation*}
      \Vert (K_c \omega)_\ell \Vert _{W_c^{1,p}} \lesssim \big \Vert (K_c \omega)_\ell \big\Vert_{p} + \big\Vert|\nabla|(K_c \omega)_\ell \big \Vert_{p}  \lesssim \big\Vert |\nabla|^{-\ell}\omega \big \Vert_p + \big \Vert |\nabla|^{1-\ell} \omega \big \Vert_p,
    \end{equation*}
    According to Proposition \ref{pro:nablaminus}, the right hand side is bounded by $\|\omega\|_{q_\ell}+\|\omega\|_{q_{\ell-1}}$, where $p^{-1}=q_s^{-1}-(s/Q)$ for $s=0,\ldots,\delta$. The smallest exponent that arises is $q_{\delta}$. By density of $\mathcal{S}_0$ in all $L^q$ (\cite[Proposition 7.1]{PansuRumin2018AHL}), $K_c$ extends to a bounded operator from $L^{p}\cap L^{q}(\mathbb{G})$ to $W_{c}^{1,p}(\mathbb{G})$. Proposition \ref{pro:sobolev-embedding} allows to replace this Sobolev space with $\dot C_{c}^{1-(Q/p)}(\mathbb{G})\cap L^{p}(\mathbb{G})$.
\end{pf}

   \subsection{Localization}
   
   We borrow the following trick from \cite{BaldiFranchiPansu2021Cohomology}. The left-invariant operator $K_c$ is given by convolution with a function $k_c$ which is smooth away from the origin in $\mathbb{G}$. Let $\chi$ be a smooth cut-off function, with support in a small ball (its diameter should be at most the difference between the radii of $B$ and $B'$) and which is equal to $1$ in a neighborhood of the origin. Let $k=\chi\,k_c$ and $k'=(1-\chi)\,k_c$, so that $k$ has small support and $k'$ is smooth on $\mathbb{G}$. Let $P$ (resp. $P'$) denote the operator of convolution with $k$ (resp. $k'$), so that $K_c=P+P'$. Then
   $$
   1=\d P+P\d + S,
   $$
where $S=\d P'+P' \d$ is smoothing. Let $\rho$ denote the operator of restriction from $B$ to $B'$. We note that $\rho P$ is well defined from $\mathscr{E}^{m}_{c}(B)$ to $\mathscr{E}^{m-1}_{c}(B')$. As a consequence, $\rho S=\rho-\d \rho P-\rho P\d$ is well defined and smoothing from $\mathscr{E}^{m}_{c}(B)$ to $\mathscr{E}^{m}_{c}(B')$.

  \begin{lem}\label{lem:sobolev-local}
Fix $p\in(1,\infty)$. If $\omega\in \mathscr{E}^{m}_{c}(B)$, then
\begin{align*}
\|\rho P\omega\|_{C_{c}^{1-(Q/p)}(B')}&\lesssim \|\omega\|_{C^0(B)},\\
\|\rho P\d\omega\|_{C_{c}^{1-(Q/p)}(B')}&\lesssim \|\d\omega\|_{C^0(B)},\\
\|\rho S\omega\|_{C_{c}^{1-(Q/p)}(B')}&\lesssim \|\omega\|_{C^0(B)}+\|\d\omega\|_{C^0(B)}
\end{align*}
In particular,
\begin{align*}
\|\rho P\omega\|_{C_{c}^{1-(Q/p)}(B')}+
\|\rho (P\d+S)\omega\|_{C_{c}^{1-(Q/p)}(B')}&\lesssim \mathbf{F}_{B}(\omega).
\end{align*}
  \end{lem}
  
  \begin{pf}
  Let us view $L^p(B)$ as a subspace of $L^p(\mathbb{G})$ by extending $L^p$ forms on $B$ by $0$. Since it has a smooth kernel, $\rho P'$ is bounded from $L^p(B)$ to any H\"older space, for instance $C_{c}^{k,1-(Q/p)}(B')$ for $k$ larger than the order of $\d$ plus one, whence the estimate 
$$
\|\rho S\omega\|_{C_{c}^{1-(Q/p)}(B')\cap L^p(B')}\lesssim \|\omega\|_{L^p(B)}+\|\d\omega\|_{L^p(B)}\lesssim \|\omega\|_{C^0(B)}+\|\d\omega\|_{C^0(B)}.
$$
Also, according to Lemma \ref{lem:kc}, $\rho P=\rho K_c -\rho P'$ is bounded from an intersection of spaces $L^q(B)\subset L^q(\mathbb{G})$ to $\dot C_{c}^{1-(Q/p)}(B')\cap L^p(B')=C_{c}^{1-(Q/p)}(B')$. This intersection contains $C^0(B)$. 
\end{pf}
   
\subsection{Proof of Theorem \ref{thm:forms}}
  
  \begin{pf}
   Consider a bounded sequence $e_j$ of classes in $E_1$. We can pick representatives $\omega_j$, and by approximation, we can choose smooth Rumin forms $\omega'_j$ defined on $B$, uniformly bounded in the flat norm $\mathbf{F}_B$, and such that $\mathbf{F}_{B}(\omega'_j-\omega_j)$ tends to $0$.
   
   According to Lemma \ref{lem:sobolev-local}, the restrictions to the smaller ball $B'$ can be written
   \begin{align*}
\rho(\omega'_j)=\rho (P\d +S)(\omega'_j)+\d\rho P(\omega'_j),
\end{align*}
where $\phi_j:=(P\d +S)(\omega'_j)$ and $\psi_j:=\rho P(\omega'_j)$ are bounded in $C_{c}^{1-(Q/p)}(B')$. By Arzel\`a-Ascoli, up to extracting a subsequence, one can assume that $\phi_j$ and $\psi_j$ converge in $C^0(B')$. Therefore the subsequence $\rho(\omega'_j)$ converges in the $\mathbf{N}_{B'}$ norm. So does a subsequence of the original sequence $\omega_j$, when restricted to $B'$. This shows that a subsequence of $\Theta(e_j)$ converges in $\mathbf{F}_{c,m}(\mathbb{G})$.

\end{pf}

\section{Charges}

The space $CH$ of Rumin charges is one of the many possible duals to the locally convex topological space of Rumin currents with compact support on $\mathbb{G}$. But it turns out to have a simple description.

If $\phi,\psi$ are continuous Rumin forms on $\mathbb{G}$, then $\phi+\d\psi$ obviously defines a charge. Conversely, we shall show that every Rumin charge is of this form.

Since the discussion of Rumin charges in Carnot groups exactly parallels that of charges in Euclidean spaces, we merely survey \cite{dePauwMoonensPfeffer2009Charges}, highlighting where Theorem \ref{thm:compactness} is needed.

\subsection{Semireflexivity}

Fix a compact set $K\subset \mathbb{G}$, let $\mathcal{M}_K$ (resp. $\mathcal{N}_K$) be the space of currents of finite mass (resp. finite $N$ mass) with support in $K$. Then $\mathcal{M}_K$ is the dual of the space of continuous Rumin forms on $K$,  $\mathcal{M}_K= C^0(K)'$. It follows that $\mathcal{N}_K$ is a dual as well,
$$
\mathcal{N}_K=\mathcal{M}_K \cap \d^{-1}(\mathcal{M}_K)=(C^0(K)+\d C^0(K))'.
$$
So we are aiming at a form of reflexivity of $C^0+\d C^0$, which fails!  The point is the change of topology on $\mathcal{N}_K$, passing from the $\mathbf{N}$ mass to the flat mass. The keyword is \emph{semireflexivity},  
\cite[Chapter IV.2]{Bourbaki}.

\begin{dfi}[Bourbaki]
A linear form on a topological vectorspace is \emph{strongly continuous} if it is bounded on bounded subsets.

A locally convex topological vectorspace $X$ is \emph{semireflexive} if every strongly continuous linear form on its topological dual $X'$ arises from an element of $X$.
\end{dfi}
%


Let $X$ be a locally convex space, let $\mathcal{S}$ be a dilation stable family of convex subsets of $X$. There is a topology $\mathcal{T}_{S}$ on $X$, inducing the initial topology on each $S\in\mathcal{S}$, such that a linear map $X\to Y$ is continuous if and only if all its restrictions to elements of $S$ are continuous, \cite[Proposition 4.3]{dePauwMoonensPfeffer2009Charges}.

\begin{pro}[{\cite[Theorem 3.16]{dePauwMoonensPfeffer2009Charges}}]
\label{pro:semireflexive}
Let $X$ be a locally convex space, let $\mathcal{S}$ be an exhausting family of \emph{compact} convex subsets of $X$. Then $(X,\mathcal{T}_{S})$ is semireflexive.
\end{pro}

If $\mathcal{N}$ is the space of compactly supported normal Rumin currents, endowed with the flat topology, and $\mathcal{S}=\{S_{K,\nu}\}$ (see the definition in Section \ref{sec:intro-results}), then Rumin charges are exactly the continuous linear functionals on $(\mathcal{N},\mathcal{T}_{S})$,
$$
CH\simeq(\mathcal{N},\mathcal{T}_{S})'.
$$

Thanks to Theorem \ref{thm:compactness}, Proposition \ref{pro:semireflexive} applies. By semireflexivity, every strongly continuous linear functional on the space of charges $CH$ arises from a normal current with compact support,
$$
(\mathcal{N},\mathcal{T}_{S})\simeq CH^*.
$$

\subsection{Proof of the representation of charges}

The second ingredient is the identification of bounded subsets of $C^0$.

\begin{pro}[{\cite[Lemma 6.4]{dePauwMoonensPfeffer2009Charges}}]
If a subset $S$ of $\mathcal{N}$ is uniformly bounded as linear functionals on $C^0$, then all elements of $S$ have support in the same compact set of $\mathbb{G}$. 
\end{pro}

We want to show that $\Theta:C^0\oplus C^0\to CH$, $(\phi,\psi)\mapsto \phi+\d\psi$ is onto. We proceed by showing that its adjoint $\Theta^*:\mathcal{N}=CH^*\to (C^0\oplus C^0)^*$, given by
$$
\langle \Theta^*(T),(\phi,\psi) \rangle = \langle T,\phi+\d\psi \rangle=\langle T,\phi \rangle+\langle \partial_c T,\psi \rangle,
$$
is proper.  If $S\subset \mathcal{N}$ and $\Theta^*(S)$ is bounded on all bounded subsets of $C^0\oplus C^0$, then, for $T\in S$, $N(T)$ is bounded and $\mathrm{supp}(T)$ is in a common compact set. By Theorem \ref{thm:compactness}, $S$ is flat-compact. 

By \cite[Proposition 6.8]{dePauwMoonensPfeffer2009Charges}, this implies that the range of $\Theta^*$ is weak$^*$-closed, and hence, by the closed range theorem, that the range of $\Theta$ is closed. This range contains the dense subspace of smooth compactly supported Rumin forms, so $\Theta$ is onto, i.e. every charge can be written $\phi+\d\psi$, for $\phi,\psi\in C^0$.

\bibliographystyle{abbrv}
\bibliography{./references}

\end{document}